\documentclass[12pt]{article}
\usepackage{amssymb}

\usepackage[dvips]{graphicx}
\usepackage{amsmath}

\newtheorem{theorem}{Theorem}[section]

\newtheorem{corollary}[theorem]{Corollary}

\newenvironment{proof}[1][Proof]{\textbf{#1.} }{\ \rule{0.5em}{0.5em}}

\begin{document}
\title{Matchings in arbitrary groups}
\date{}
\author{Shalom Eliahou and C\'{e}dric Lecouvey \\
Laboratoire de Math\'{e}matiques Pures et Appliqu\'{e}es \\
Universit\'{e} du Littoral C\^{o}te d'Opale\\
50 rue F. Buisson, B.P. 699\\
62228 Calais Cedex, France}
\maketitle

\begin{abstract}
A \textit{matching} in a group $G$ is a bijection
$\varphi$ from a subset $A$ to a subset $B$ in $G$ 
such that $a \varphi(a) \notin A$ for all $a \in A$. 
The group $G$ is said to have the \textit{matching property} if,
for any finite subsets $A,B$ in $G$ of same
cardinality with $1 \notin B$, there is a matching from
$A$ to $B$. 

Using tools from additive number theory,  Losonczy
proved a few years ago 
that the only abelian groups satisfying the matching property 
are the torsion-free ones and those of prime order. He also proved that, 
in an abelian group, any finite subset $A$ avoiding $1$ admits a matching 
from $A$ to $A$. 

In this paper, we show that both Losonczy's results hold verbatim 
for arbitrary groups, not only abelian ones. Our main tools
are classical theorems of Kemperman and Olson, also
pertaining to additive number theory, but  
specifically developped for possibly nonabelian groups.
\end{abstract}

\section{Introduction}

Let $G$ be a group, written multiplicatively. Given nonempty finite subsets $A,B$ in $G$, 
a \textit{matching} from $A$ to $B$ is a map $\varphi:A \rightarrow B$
which is bijective and satisfies the condition $$a\varphi(a) \notin A$$ for all $a \in A$. 

This notion was introduced in \cite{FanLos} by Fan and Losonczy, who used matchings in $\mathbb{Z}^n$ 
as a tool for studying an old problem of Wakeford concerning canonical forms for symmetric
tensors \cite{Wak}.

Coming back to general groups, it is plain that if there is a matching $\varphi$ from $A$ to $B$, then $|A|=|B|$ and $1 \notin B$. (For if $1 \in B$, let $a_1 = \varphi^{-1}(1)$; then  $a_1 \varphi(a_1)= a_1 \in A$.) It is natural to wonder whether these necessary conditions for the existence of a matching from $A$ to $B$ are also sufficient. The answer turns out to depend on the group structure. 

Following Losonczy, we say that the group $G$ has the \textit{matching property} if, 
whenever the subsets $A,B$ satisfy the conditions $|A|=|B|$ and $1 \notin B$, there exists
a matching from $A$ to $B$. Losonczy proved the following result.

\begin{theorem}[\cite{LOs}]\label{thm:matching property} Let $G$ be an abelian group. Then $G$ has the matching property if and only if $G$ is torsion-free or cyclic of prime order.
\end{theorem}

A special case of interest is the one where $A=B$. Is it sufficient, in this case, to assume that $A$ does not contain $1$ in order to guarantee the existence of a matching from $A$ to $A$? Losonczy's answer for abelian groups is yes.

\begin{theorem}[\cite{LOs}]\label{thm:automatching} Let $G$ be an abelian group. Let $A$ be a nonempty finite subset of $G$. Then there is a matching from $A$ to $A$ if and only if $1 \notin A$.
\end{theorem} 

The proofs in \cite{LOs} are based on methods and results from additive number theory, namely the Dyson transform, and theorems of Cauchy-Davenport and Kneser. However powerful, these methods only work for abelian groups.

\medskip

In Section~\ref{secM} of this paper, we extend the above two theorems of Losonczy to arbitrary groups. This is achieved by making use of results in additive number theory which were specifically developped for possibly nonabelian groups. These results are recalled in the next section. The engine behind their proofs is the Kemperman transform, a clever nonabelian analogue of the Dyson transform. See Olson's paper \cite{Ols}.
See also Nathanson's book \cite{Nath} for general background on additive number theory.

\section{Nonabelian additive theory}

Given subsets $A,B$ of a group $G$, their \textit{product set} is defined as $$AB=\{ab \mid a \in A, b \in B\}.$$ We start with a result of Kemperman providing a conditional lower bound on the size of $AB$. 

\begin{theorem}[Kemperman \cite{Kem}]
\label{th_ke}\label{th_kem}\label{THK}
Let $A,B$ be finite subsets of a group $G.$ Assume there
exists an element $c\in AB$ appearing exactly once as a product $c=ab$ with $a\in A, \, b \in B$.
Then 
\begin{equation*}
\left| AB\right| \geq \left| A\right| +\left| B\right| -1.
\end{equation*}
\end{theorem}

The following corollary will be used in the next section for our extension of Theorem~\ref{thm:automatching}.

\begin{corollary}
\label{corKem} Let $U,V$ be nonempty finite subsets of a group $G$ such that
$U$, $V$ and $UV$ are all three contained in a subset $X$ of $G \setminus \{1\}.$ Then
\begin{equation*}
\left| X \right| \ge \left| U \right| + \left| V \right| + 1.
\end{equation*}
\end{corollary}
\begin{proof}
Let $A = U\cup \{1\}$, $B=V\cup \{1\}$.\ Then $1 \in AB$
and appears exactly once as a product in $AB$. Indeed, 
assume $1=ab$ with $a \in A$, $b \in B$. Then either $a=1$ or $b=1$, 
since $1 \notin UV$ by hypothesis, and hence $a=b=1$.
Therefore Theorem~\ref{THK} applies, and gives
$$
\left| AB \right| \ge \left| A \right| + \left| B \right| -1.
$$
Since $|A|=|U|+1$, $|B|=|V|+1$ and $AB=UV \cup U \cup V$, we have $AB \subset X$ and hence
$$
\left| X \right| \ge \left| AB \right| \ge \left| U \right| + \left| V \right| +1,
$$ 
as desired.
\end{proof}

\medskip

As for extending Theorem~\ref{thm:matching property} to arbitrary groups, we shall need
the following result of Olson.

\begin{theorem}[Olson \cite{Ols}]\label{thOl}
Let $A,B$ be nonempty finite subsets of a group $G$. There exists a finite subgroup $H$ of $G$
and a nonempty subset $T$ of $AB$ such that 
\begin{equation*}
\left| AB\right| \ge \left| T\right| \geq \left| A\right| +\left| B\right| -\left| H\right|, 
\end{equation*}
and either $HT=T$ or $TH=T$.
\end{theorem}

\section{Results and proofs\label{secM}}

We now present our extensions of Losonczy's theorems. 
Besides the additive tools from the preceding section, we shall also need, as in \cite{FanLos, LOs}, 
the marriage theorem of Hall. Recall that,
given a collection $\mathcal{E}=\{E_{1},E_{2},\dots,E_{n}\}$ of subsets of a set $E$, 
a \textit{system of distinct representatives} for $\mathcal{E}$ is a set 
$\{x_{1},\dots,x_{n}\}$ of pairwise distinct elements of $E$ with the
property that $x_{i}\in E_{i}$ for all $i=1,\dots,n.$ Hall's theorem gives necessary
and sufficient conditions for the existence of such systems. 

\begin{theorem}[Hall \cite{Hall}]
Let $E$ be a set and $\mathcal{E}=\{E_{1},E_{2},\dots,E_{n}\}$ 
a family of finite subsets of $E.$ Then $\mathcal{E}$ admits a system of distinct
representatives if and only if 
\begin{equation*}
\left| \bigcup_{i \in S}E_{i}\right| \geq \left| S\right|
\end{equation*}
for all nonempty subsets $S \subset \{1,\dots,n\}.$\end{theorem}

We are now ready to generalize Theorem~\ref{thm:automatching}.

\begin{theorem}
\label{th-match1}Let $G$ be a group. Let $A$ be a nonempty finite subset of $G$. Then there is
a matching from $A$ to $A$ if and only if $1 \notin A$.
\end{theorem}
\begin{proof} We already know that if $A$ contains $1$, there cannot be a matching from $A$ to $A$. 
Assume now $1 \notin A$. For each $a\in A$, set 
\begin{equation*}
E_{a}=\{x \in A\mid ax \notin A\}.
\end{equation*}
Finding a matching from $A$ to $A$ is clearly equivalent to finding a system of
distinct representatives for the family of sets
$$\mathcal{E} = \{E_a \mid a \in A \}.$$
By the Hall marriage theorem, this is also equivalent to the inequalities
\begin{equation}
\left| \bigcup_{s\in S}E_{s}\right| \geq \left| S\right|   \label{lneqHall}
\end{equation}
for all nonempty subsets $S\subset A$. 

Denote $E_s' = A \setminus E_s$, 
the complement of $E_s$ in $A$. Hall's conditions (\ref{lneqHall})
may be rewritten as
\begin{equation}\label{Hall2}
\left| \bigcap_{s\in S}E_s'\right| \leq \left|A\right|-\left| S\right|
\end{equation}
for all nonempty subsets $S\subset A$. Set 
$$V_S=\bigcap_{s\in S}E_s' = \{x\in A\mid sx\in A \textrm{ for all } s\in S\}.$$
We have $S V_S \subset A$ by construction. Since $1\notin A$, Corollary \ref{corKem} applies
(with $U,V,X$ standing for $S,V_S, A$ respectively), and gives 
$$\left| S\right| +\left| V_S\right| \le \left| A\right|-1 .$$
This shows that conditions (\ref{Hall2}) are satisfied and finishes the proof of the theorem.
\end{proof}

\bigskip

We now turn to the characterization of all groups satisfying the matching property.
The abelian case was first settled by Losonczy as Theorem~\ref{thm:matching property}.

\begin{theorem} Let $G$ be any group. Then $G$ has the matching property
if and only if $G$ is torsion-free or cyclic of prime order.
\end{theorem}
\begin{proof}
Assume first that $G$ is neither torsion-free nor cyclic of prime order. 
Then there is an element $a \in G$, of finite order $n \ge 2$, which does not generate $G$. 
Let $$A = \langle a \rangle = \{1,a,\dots,a^{n-1}\}$$
be the subgroup generated by $a$. Let $g \in G \setminus A$ and set 
$$
B = A \cup \{g\} \setminus \{1\} = \{a,\dots,a^{n-1},g\}.
$$
Let $\varphi: A \rightarrow B$ be any bijection. Can it possibly satisfy
the condition $x\varphi(x) \notin A$ for all $x \in A$? No, it cannot.
Picking $a \in B$ and $x_0=\varphi^{-1}(a) \in A$, we have 
$x_0 \varphi(x_0) = x_0 a \in A$ since $A$ is a subgroup. We conclude
that $G$ does not satisfy the matching property.

\smallskip

Conversely, assume that $G$ is either torsion-free or cyclic of prime order.
This means that the only finite subgroups of $G$ are $\{1\}$, and $G$ if $G$ is finite.
The trivial group is torsion-free and vacuously satisfies the matching property.
Assume now $G \neq \{1\}$.
Let $A,B$ be nonempty finite subsets of $G$ with $|A|=|B|$ and $1 \notin B$. For each $a\in A$, set 
\begin{equation*}
E_{a}=\{x \in B\mid ax \notin A\}.
\end{equation*}
Again, finding a matching from $A$ to $B$ is equivalent to finding a system of
distinct representatives for the family of sets 
$$\mathcal{E} = \{E_a \mid a \in A \}.$$
By the Hall marriage theorem, it suffices to prove the inequalities
\begin{equation}
\left| \bigcup_{s\in S}E_{s}\right| \geq \left| S\right|   \label{lneqHall3}
\end{equation}
for all nonempty subsets $S\subset A$. Denote $E_s' = B \setminus E_s$, 
the complement of $E_s$ in $B$. Hall's conditions (\ref{lneqHall3})
may be rewritten as
\begin{equation}\label{Hall4}
\left| \bigcap_{s\in S}E_s'\right| \leq \left|A\right|-\left| S\right|
\end{equation}
for all nonempty subsets $S\subset A$. Set 
$$V_S=\bigcap_{s\in S}E_s' = \{x\in B\mid sx\in A \textrm{ for all } s\in S\},$$
and $W_S = V_S \cup \{1\}$. We have $|W_S|=|V_S|+1$ and $S W_S \subset A$ by construction. 
By Theorem \ref{thOl}, there is a finite subgroup $H \subset G$ and a nonempty subset $T \subset S W_S$ such that
\begin{equation}\label{Hall5}
|S W_S| \geq |S| + |W_S| - |H|
\end{equation}
and $HT=T$ or $TH=T$. We cannot have $H=G$, for otherwise $T=G$. But as
$T \subset S W_S \subset A$, this would imply $A=G=B$, contradicting the hypothesis $1 \notin B$.
It follows that $H = \{1\}$, and inequality (\ref{Hall5}) yields
$$
|A| \geq |S| + |V_S|,
$$
since $S W_S \subset A$, $|W_S|=|V_S|+1$ and $|H|=1$.
Therefore conditions (\ref{Hall4}), which imply the existence of a matching
from $A$ to $B$, are satisfied. It follows that $G$ has the matching property.
\end{proof}

\end{document}